\documentclass[times]{nlaauth}

\usepackage{moreverb}

\usepackage[dvips,colorlinks,citecolor=red,urlcolor=red]{hyperref}

\def\volumeyear{2011}

\begin{document}

\def\RR {\mathbb{R}}
\def\NN {\mathbb{N}}
\def\CC {\mathbb{C}}
\def\KK {\mathcal{K}}
\def\II {\mathcal{I}}
\def\vep {\varepsilon}
\def\UU {\mathcal{U}}
\def\VV {\mathcal{V}}

\newcommand{\hdrie}{\hspace{3mm}}
\newcommand{\und}{\hdrie\mbox{\rm and }\hdrie}
\newcommand{\half}{\frac{1}{2}}
\newcommand{\be}{\begin{equation}}
\newcommand{\ee}{\end{equation}}
\newcommand{\sth}{\hdrie | \hdrie}
\newcommand{\rank}{{\rm rank}}
\newcommand{\sspan}{{\rm span}}

\newtheorem{Theorem}{Theorem}[section]
\newtheorem{Lemma}[Theorem]{Lemma}
\newtheorem{Corollary}[Theorem]{Corollary}
\newtheorem{Proposition}[Theorem]{Proposition}
\newtheorem{Definition}[Theorem]{Definition}
\newtheorem{Remark}[Theorem]{Remark}

\runningheads{J.~H.~Brandts and R.~Reis da Silva}{On the SPAM method}

\title{On the Subspace Projected Approximate Matrix method}

\author{J.~H.~Brandts\corrauth and R.~Reis da Silva}

\address{Korteweg-de Vries Institute for Mathematics, Faculty of Science, University of Amsterdam, Netherlands}

\corraddr{Korteweg-de Vries Institute for Mathematics, Faculty of Science, University of Amsterdam, P.O. Box 94248, 1090 GE Amsterdam, Netherlands. E-mail: J.H.Brandts@uva.nl}

\begin{abstract}
We provide a comparative study of the Subspace Projected Approximate Matrix method, abbreviated SPAM, which is a fairly recent iterative method to compute a few eigenvalues of a Hermitian matrix $A$. It falls in the category of inner-outer iteration methods and aims to save on the costs of matrix-vector products with $A$ within its inner iteration. This is done by choosing an approximation $A_0$ of $A$, and then, based on both $A$ and $A_0$, to define a sequence $(A_k)_{k=0}^n$ of matrices that increasingly better approximate $A$ as the process progresses. Then the matrix $A_k$ is used in the $k$th inner iteration instead of $A$.

In spite of its main idea being refreshingly new and interesting, SPAM has not yet been studied in detail by the numerical linear algebra community. We would like to change this by explaining the method, and to show that for certain special choices for $A_0$, SPAM turns out to be mathematically equivalent to known eigenvalue methods. More sophisticated approximations $A_0$ turn SPAM into a boosted version of Lanczos, whereas it can also be interpreted as an attempt to enhance a certain instance of the preconditioned Jacobi-Davidson method.

Numerical experiments are performed that are specifically tailored to illustrate certain aspects of SPAM and its variations. For experiments that test the practical performance of SPAM in comparison with other methods, we refer to other sources. The main conclusion is that SPAM provides a natural transition between the Lanczos method and one-step preconditioned Jacobi-Davidson.
\end{abstract}


\maketitle

\section{Introduction}

We provide a comparative study of SPAM \cite{ShWaTiMi}. SPAM, which stands for Subspace Projected Approximate Matrix, is a fairly recent (2001) method for the computation of eigenvalues of a large Hermitian matrix $A$. Like the Davidson method \cite{Dav}, SPAM was originally developed for matrices that arise from applications in Chemistry. It was only in \cite{CrPhSa}, many years after its conception, that Davidson's method was given proper attention by the numerical linear algebra community. As far as we can tell, also SPAM has been neglected by the numerical linear algebra community. Moreover, even though a number of citations \cite{ChPo,GySeCh,MeGrWaMiSh,RiLuLe,ZhShMi} within the Chemistry and Physics communities over the past years demonstrate awareness of its existence, no studies of its mathematical properties seems to exists.

SPAM belongs to the category of inner-outer iteration methods and is interpreted by its inventors as a modification of the above mentioned Davidson method. It is based on the following observation. Even when sparse, the computational effort necessary to carry out a matrix-vector multiplication with $A$ can be significant and often represents the bottleneck of the total computational effort. SPAM reduces the costs of matrix-vector products with $A$ by replacing its action within the inner iteration of the algorithm with a sparser or more structured approximation. By doing so, it attempts to slash the overall computational cost. The idea is not altogether new and is related to a certain type of preconditioning, called {\em one-step approximation} in the Jacobi-Davidson method. See Section 4.1 of \cite{SlVo}. There too, the matrix $A$ is, in the inner iteration, replaced by a preconditioner. The originality of the approach in \cite{ShWaTiMi} lies in the fact that the action of the preconditioner is only applied to the subspace in which the action of $A$ has not yet been computed in the outer iteration of the method. Consequently, the approximated action of $A$ in the inner iterations is likely to become more and more accurate as the number of outer iterations increases. Note that nonetheless, only one approximate matrix $A_0$ is needed. Intuitively, one would expect that SPAM would outperform Jacobi-Davidson with one-step approximation.

Since the main idea of SPAM is refreshingly new and potentially interesting, with links to other eigenvalue methods and selection techniques, we would like to bring it to the attention of the numerical linear algebra community. Indeed, following an initial studies that led to a MSc thesis by the second author, there is a need to understand the method in more detail, not only in theory, but also to submit it to a more illustrative set of numerical experiments than those in \cite{ShWaTiMi}. In the experiments in \cite{ShWaTiMi}, the overall efficiency of SPAM was the main interest instead of its mathematical position within the class of iterative eigensolvers. In particular, we would like to point out the similarities with and differences to strongly related and well-known iterative methods such as the Lanczos method \cite{Lan} and the Jacobi-Davidson method of Sleijpen and van der Vorst \cite{SlVo}.

\subsection{Outline of the results in this paper}

We will show that for certain choices of the approximate matrix $A_0$, the SPAM method is mathematically equivalent to methods such as Lanczos \cite{Lan} or the Riccati \cite{Bra} method, another attempt to improve upon the Jacobi-Davidson \cite{SlVo} method. Further, we will see that a Schur complement-based choice for the action of $A$ outside the Ritz-Galerkin subspace that is being built in the outer iteration naturally leads to a connection with harmonic Rayleigh-Ritz \cite{Bea} methods. Next, we show that choosing $A_0$ such that $A-A_0$ is positive semi-definite has, at least in theory, an agreeable effect on the approximations obtained in the inner iteration in comparison to choosing $A_0$ without such a restriction. Numerical experiments suggest that this also works through into the outer iteration. We comment on how such approximations $A_0$ can be obtained in the context of discretized elliptic PDEs \cite{CiLi,WrAl} but also from a purely algebraic point of view. Finally, we present a variety of detailed numerical illustrations of the performance of the method in comparison with the Lanczos and the Jacobi-Davidson method. These illustrations do not merely aim to show that SPAM is a suitable method to solve eigenvalue problems (which was for a large part already taken care of in \cite{ShWaTiMi}), but to emphasize the role of approximations $A_0$ from below and to show the similarities and discrepancies with Lanczos and Jacobi-Davidson, such that the SPAM method can be put into a proper perspective.

\section{SPAM and some other subspace methods}

Eigenvalue problems are among the most prolific topics in Numerical Linear Algebra \cite{BaEtc,GoVo}. In particular, the continuous increase in matrix sizes turns the understanding of known iterative methods, as well as the development of more efficient ones, into an ample field of work within the global topic. Successful methods like the Implicitly Restarted Arnoldi \cite{Sor} and the Krylov-Schur \cite{Ste} methods and their symmetric counterparts are nowadays among the most competitive. For convenience and in order to set the ground for what is to come, we opt to outline the main concepts. For more detailed considerations on both theoretical and practical aspects of the numerical solution of large eigenvalue problems, we refer to \cite{Par,GoLo,Ste2,StSu}.

\subsection{Ritz values and vectors}

Throughout this paper, $A$ is a Hermitian $n\times n$ matrix with eigenvalues
\be \label{eigs} \lambda_n \leq \lambda_{n-1} \leq \cdots \leq \lambda_2 \leq \lambda_1, \ee
Let $V$ be an $n\times k$ matrix with mutually orthonormal columns and define the $k$-dimensional subspace $\VV$ of $\CC^n$ by
\be \VV = \{ Vy \sth y\in \CC^k\}. \ee
In the context of iterative methods, $\mathcal{V}$ is called the \emph{search subspace}. Let $V_{\perp}$ be such that $(V|V_{\perp})$ is unitary and write
\be\label{eq:Asubsp}
\hat{A} =   (V|V_{\perp})^*A(V|V_{\perp}) = \left[\begin{array}{cl}
    M & R^{*} \\
    R & S
\end{array}\right].
\ee
The eigenvalues of the $k\times k$ matrix $M=V^*AV$,
\be \label{ritz} \mu_k \leq \mu_{k-1} \leq \cdots \leq \mu_2 \leq \mu_1, \ee
are called the \emph{Ritz values} of $A$ with respect to $\mathcal{V}$. The vectors $u_{i}=Vz_{i}$, where $z_1,\dots,z_k$ is an orthonormal basis for $\CC^k$ consisting of eigenvectors of $M$ belonging to the corresponding $\mu_i$, are the \emph{Ritz vectors} of $A$ in $\VV$. The \emph{residuals} $\hat{r}_i=Au_i-u_i\mu_i$ for the respective Ritz pairs $(\mu_i,u_i)$ satisfy
\be \label{two} Au_i-u_i\mu_i = \hat{r}_i\perp\VV.\ee
Each \emph{Ritz pair} $(\mu_i,u_i)$ is also an eigenpair of the $n\times n$ rank-$k$ matrix $VMV^*$ and is interpreted as an approximation of an eigenpair of $A$. See \cite{JiSt,Par,SlEs}.

\subsection{Rayleigh-Ritz and subspace expansion}

The \emph{Rayleigh-Ritz} procedure is the first stage in iterative methods for eigenproblems and consists of computing the Ritz pairs from $\VV$.  The computation of $S$ in (\ref{eq:Asubsp}) is not needed, nor feasible for reasons of efficiency. However, a cheaply available by-product of the computation of $AV$ is the matrix $\hat{R}=AV-VM$, where $R=V_\perp^*\hat{R}$. Its columns are the respective residuals (\ref{two}). In the second stage, the search subspace $\VV$ is expanded. Different definitions of the expansion vector distinguish the different iterative methods. Each strategy results in a sequence of nested spaces
\be \VV_1 \subset \VV_2 \subset \dots \subset \VV_{n-1} \subset \VV_n\ee
and has the objective to obtain accurate Ritz pairs while spending only minimal computational effort. One of the strategies will result in SPAM. Other strategies lead to methods with which SPAM will be compared in this paper: the Lanczos \cite{Lan}, Jacobi-Davidson \cite{SlVo}, and Riccati \cite{Bra} methods.

\subsection{The Lanczos method}\label{sect-Lanczos}

The Lanczos method \cite{Lan} defines $\VV_{j+1}$ as $\VV_j\oplus\sspan\{r\}$ where $r\perp\VV_j$ is any of the current residuals from (\ref{two}). This results in a well-defined method: starting with some initial vector $v_1$ with $\|v_1\|=1$ that spans the one-dimensional search space $\VV_1$, it can be easily verified by induction that regardless which residual is used for expansion, the sequence of search spaces that is defined, equals the sequence of Krylov subspaces
\be \VV_j = \KK^j(A,v_1) = \sspan\{v_1, Av_1, \dots, A^{j-1}v_1\}. \ee
Due to (\ref{two}), the matrix $V_k$ with the property that its column span equals $\KK^k(A,v_1)$ can be chosen to have as columns $v_1$ and the normalized residuals
\be v_{j+1} = \frac{\hat{r}_j}{\|\hat{r}_j\|}, \hdrie j\in\{1,\dots,k\},\ee
where $\hat{r}_j$ is a residual from the $j$-dimensional search space. From this the so-called \emph{Lanczos relation} results,
\be\label{vijf} AV_k = V_{k+1}M_{k,k+1}, \hdrie\mbox{\rm with } \hdrie M_{k,k+1} = \left[\begin{array}{c} M_k \\ \hline \rho_k e_k^*\end{array}\right].\ee
Here, $M_k$ is a $k\times k$ tridiagonal matrix, $e_k$ is the $k$th canonical basis vector of $\CC^k$ and $\rho_k$ is a scalar. The following trivial observation will have its counterpart in the discussion of the SPAM method in Section \ref{sectSPAM}.

\begin{Remark}\label{remark1} If the Lanczos method runs for the full number of $n$ iterations, it produces a unitary matrix $V_n$ that depends only on $A$ and the start vector $\|v_1\|$. The $k$th leading principal submatrix of the $n\times n$ tridiagonal matrix $M=V_n^*AV_n$ is the matrix $M_k$ from (\ref{vijf}) whose eigenvalues are the Ritz values after $k$ iteration steps.
\end{Remark}

\subsection{The Subspace Projected Approximate Matrix (SPAM) method}\label{sectSPAM}

In the Subspace Projected Approximate Matrix (SPAM) method \cite{ShWaTiMi}, the expansion vector is a suitable eigenvector of an approximation of $A$. This approximation has a cheaper action than $A$ itself. Thus, the matrix-vector products within the inner iteration that is needed to compute this eigenvector, will be cheaper than for instance in the Jacobi-Davidson \cite{SlVo} and Ricatti \cite{Bra} methods. These methods, which we will explain in more detail in Section \ref{Sect2.3}, both use $A$ itself in their inner iteration.

A central observation in \cite{ShWaTiMi} is that the action of $A$ on $\VV$ that has already been performed in the outer iteration, can be stored in a matrix $W=AV$, and be re-used within the inner iteration at relatively low costs. Thus, the action of $A$ in the inner iteration only needs to be approximated partially. The resulting approximation is then different after each outer iteration step, even though only one approximate matrix $A_0$ is provided. Its action is merely used on less and less of the total space. As such, SPAM may be interpreted as a discrete homotopy method.

\subsubsection{General description of SPAM.}\label{sect2.4.1}

Let $A_0$ be an approximation of $A$. In Section \ref{sect3} we will comment on how this approximation can be chosen. For now, assume that $A_0$ is available and define, in view of (\ref{eq:Asubsp})
\be\label{acht} \hat{S} = {V_\perp^{\phantom{*}}} A_0V_\perp^*\ee
and
\be \label{eq:ASPAM} A_k = (V|V_{\perp})\hat{A}_k(V|V_{\perp})^*, \hdrie\mbox{\rm where }\hdrie  \hat{A}_k =\left[\begin{array}{cl}
    M & R^{*} \\
    R & \hat{S}
\end{array}\right]. \ee
The subscript $k$ of $A_k$ refers to the number of columns of $V$. In \cite{ShWaTiMi}, the matrix $A_k$ is called a \emph{subspace projected approximate matrix}. This is motivated by the fact that $A_kV=AV$ and $V^*A_k=V^*A$. In particular, since $M=V^*AV = V^*A_kV$, both $A_k$ and $A$ have the same Ritz pairs in $\VV$. This will be exploited to derive bounds for the eigenvalues of $A_k$ in Section \ref{sect3}. This is of interest since the search space $\VV$ in the outer iteration is expanded with an eigenvector of $A_k$. Note that the action of $A_k$ on $\VV^{\perp}$ does \emph{not} equal the action of $A_0$ on $\VV^{\perp}$. Since $A^*=A$, the action of $A_k$ on $\VV^\perp$ equals, in fact, the action of $A$ in $k$ linearly independent functionals on $\CC^n$.

With the convention that $\Pi_0=I$,  write
\be\label{elf} \Pi_k = V_{\perp}^{\phantom{*}}V_{\perp}^{*} = I- VV^{*}. \ee
This  shows, together with equations (\ref{acht}) and (\ref{eq:Asubsp}), that
\be\label{stbasis} A_k=-VV^{*}AVV^{*}+AVV^{*}+VV^{*}A+\Pi_k A_0 \Pi_k.\ee
Thus, the action of $A_k$ can benefit from the stored action of $A$ on $\VV$ as follows. With $W=AV$ we have that
\be\label{eq:MVspam} A_kv = -VMV^{*}v+WV^{*}v+VW^{*}v+\Pi_k A_0\Pi_k v \ee
where we have used (\ref{elf}) to avoid the numerically infeasible formation of $V_{\perp}$. Note that if $v\in\VV^\perp$, the first two terms vanish. In view of Remark \ref{remark1} we now observe the following.

\begin{Remark} If SPAM runs for the full $n$ iterations, it produces a unitary matrix $U_n$ that depends only on $A$ and $A_0$. The $k$th leading principal submatrix of the $n\times n$ matrix $M=U_n^*AU_n$ then contains the Rayleigh-Ritz approximations after $k$ steps of the outer iteration.
\end{Remark}

For theoretical purposes and without loss of generality, we assume that $V_\perp$ in (\ref{eq:ASPAM}) contains precisely the basis for the orthogonal complement of $\VV$ that SPAM is about to produce in future iterations. With respect to this basis, $\hat{A}_k$ is an update of $\hat{A}_{k-1}$ of {\em arrowhead type}, in the sense that
\be\label{veertien} \hat{A}_{k} - \hat{A}_{k-1} =
\left[\begin{array}{c|c|c}
 0^{\phantom*} & 0 & 0 \\
 \hline
 0^* & \tau & t^* \\
 \hline
 0^* & t & 0
 \end{array}\right]
 =
 \left(\begin{array}{c} 0 \\ \tau \\ t\end{array}\right) e_k^* + e_k \left(\begin{array}{ccc} 0^* & \tau & t^*\end{array}\right) - e_k\tau e_k^*,\ee
where each entry in the arrowhead formed by $t\in\CC^{n-k},\tau\in\RR$ and $t^*$, is the difference between the corresponding entries of $\hat{A}$ from (\ref{eq:Asubsp}) and $\hat{A}_0=(V|V_{\perp})^*A_0(V|V_{\perp})$ from (\ref{eq:ASPAM}). Thus, with respect to the basis defined by the columns of $(V|V_\perp)$, the matrix $\hat{A}_0$ simply transforms step by step into $\hat{A}$ in the sense that  after $k$ steps, the first $k$ columns and rows have changed into those of $\hat{A}$, and the resulting matrix is called $\hat{A}_k$. This is visualized in Figure 1.

\begin{figure}
\centering
\includegraphics[width=3.3cm]{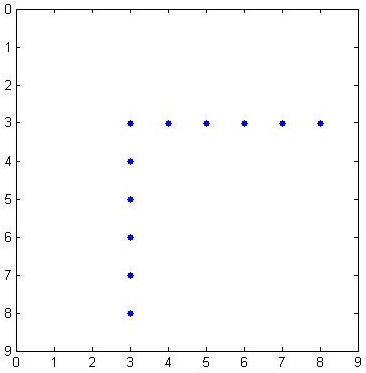} \includegraphics[width=3.3cm]{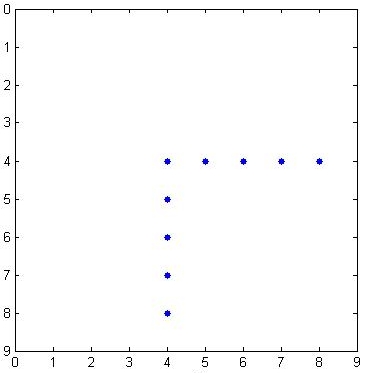} \includegraphics[width=3.3cm]{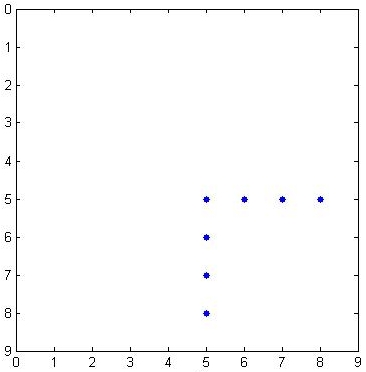} \includegraphics[width=3.3cm]{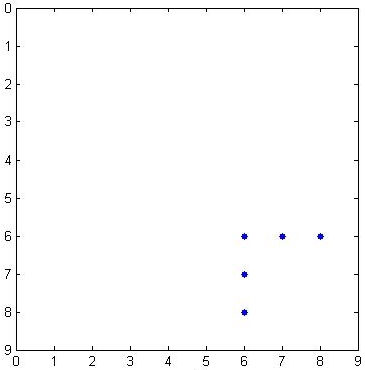}\\[2mm]
\caption{Arrowhead updates from $\hat{A}_{k-1}$ to $\hat{A}_k$ for consecutive values of $k$.}
\end{figure}

On the original basis, this transformation is described in the next proposition. Note that in this proposition,  $v$ is the eigenvector of interest of $A_{k-1}$, orthonormalized to $V_{k-1}$.

\begin{Proposition}\label{props1} Let $k\geq 1$. Write $(V_{k-1}|v|V_\perp)=(V|V_\perp)$, thus $v=Ve_k$. Then $A_k$ is the following indefinite Hermitian rank-$2$ update of $A_{k-1}$,
\be\label{update} A_k = A_{k-1} + uv^*+vu^* = A_{k-1} + (u|v)\left[\begin{array}{cc}0 & 1 \\ 1 & 0 \end{array}\right] (u|v)^*,\ee
where
\be\label{seventeenb} u = \left(\Pi_{k-1}-\half vv^*\right)(A-A_0)v =  \left(\Pi_{k}+\half vv^*\right)(A-A_0)v. \ee
\end{Proposition}

{\em Proof. } Combining (\ref{veertien}) with $(V|V_{\perp})(\hat{A}-\hat{A}_0)(V|V_{\perp})^*=A-A_0$ we find
\be \label{sixteenb}  \left(\begin{array}{c} 0 \\ \tau \\ t\end{array}\right)  = \left(\begin{array}{cccc} 0 & e_k & \cdots & e_n\end{array}\right)^*(\hat{A}-\hat{A}_0)e_k  = (0|v|V_{\perp})^*(A-A_0)v.\ee
Therefore, substituting (\ref{sixteenb}) into
\[ A_k - A_{k-1} = (V_{k-1}|v|V_{\perp}) \left( \left(\begin{array}{c} 0 \\ \tau \\ t\end{array}\right) e_k^* + e_k \left(\begin{array}{ccc} 0^* & \tau & t^*\end{array}\right) - e_k\tau e_k^*\right)(V_{k-1}|v|V_{\perp})^* \]
we arrive, using that $v=(V_{k-1}|v|V_\perp)$, at
\be A_k - A_{k-1} = \Pi_{k-1}(A-A_0)vv^* + vv^*(A-A_0)\Pi_{k-1} - vv^*(A-A_0)vv^*. \ee
Splitting the most right term in two equal parts and rearranging the terms proves the formula for $u$ in (\ref{seventeenb}). Since by (\ref{elf})
\be \Pi_{k-1}-vv^*=\Pi_k,\ee
also the second equality is proved. \hfill $\Box$

\begin{Remark} Note that the result of Proposition \ref{props1} can also be derived in an alternative way, starting from (\ref{stbasis}). The proof presented here also shows the validity of (\ref{veertien}).
\end{Remark}

From (\ref{eq:Asubsp}) and (\ref{eq:ASPAM}) we see that
\be \rank(A-A_k) = \rank(\hat{A}-\hat{A_k}) \leq n-k, \ee
and thus, even though $A_k-A_{k-1}$ has rank-$2$, the update - or maybe better downdate
\be A-A_{k}=A-A_{k-1}+(A_{k-1}-A_{k})\ee
will generically decrease the rank of $A-A_{k-1}$ with at most one. This remark goes hand in hand with the observation that even though the approximations $A_k$ of $A$ may seem unusual, the viewpoint of considering the reverse sequence $0=A-A_n, A-A_{n-1},\dots, A-A_1,A-A_0$ as increasingly better approximations of $A-A_0$ is very natural indeed: they form a sequence of Rayleigh-Ritz approximations to $A-A_0$ in the orthogonal complements of the spaces $V_k$.

\bigskip
In the outer iteration, the products $Av$ and $V^*Av$ were computed. Thus, both $v^*Av$ and $\Pi_{k-1}Av$ are available in (\ref{seventeenb}) without additional computational costs. Furthermore, since $v$ is orthogonal to the $(k-1)$-dimensional search space, we have that $v^*A_{k-1}v = v^*A_0v$. Now, because $v$ is the result of orthogonalization of an eigenvector of $A_{k-1}$ to $V_{k-1}$, also $v^*A_0v$ can be retrieved from the inner iteration. Thus, the vectors $u$ and $v$ in the updating procedure (\ref{update}) are cheaply available.

\begin{Remark} Of course, the update (\ref{update}) itself should not be performed explicitly because it will generally result in fill-in of originally sparse matrices.
\end{Remark}

\subsubsection{Choice of the method for the inner iteration of SPAM.}\label{sect2.4.2}

In \cite{ShWaTiMi}, the authors suggest to use Davidson's method \cite{Dav} to solve the eigenvalue problem for $A_k$, but of course any other method can be adopted. Apart from the Lanczos method, also the Generalized Davidson method from \cite{Mor} was tested as inner method in \cite{ShWaTiMi}. Other possibilities include, for instance, the Jacobi-Davidson \cite{SlVo} method\footnote{M. Hochstenbach presented this option at the 2006 GAMM/SIAM Applied Linear Algebra conference.}. The latter can be a good option because it often needs only a few iterations to converge if a good start vector is available. This start vector may be either the eigenvector approximation of $A_{k-1}$, or the current eigenvector approximation of $A$ from the outer iteration. We will study this choice in Section \ref{Sect2.3}.

\begin{Remark}\label{rem2.7} SPAM should first of all perform well under the assumption that the eigenproblem for $A_k$ is solved exactly. This will be investigated in the numerical illustrations in Section \ref{sect4}.
\end{Remark}

In \cite{ShWaTiMi} it is also noted that SPAM itself can be chosen in the inner iteration. This leads to a recursive multilevel version of the method, and assumes that a whole sequence of approximating matrices of $A$ is available, each having a cheaper action than its predecessor. The eigenvector computed at a given level is then used as expansion vector at the first higher level.

\subsection{Comparing SPAM with the Jacobi-Davidson and the Riccati method}\label{Sect2.3}

The philosophy of SPAM is to use all available (computed) information from (\ref{eq:Asubsp}), i.e., $M,R$ and $R^*$, to determine the expansion vector, and consequently, that only $S$ needs to be approximated. The Jacobi-Davidson \cite{SlVo} and Riccati \cite{Bra} methods partially share this philosophy. Instead of $R$, they only use the residual corresponding to the selected eigenvector approximation. On the other hand, in their most simple forms, they do not approximate the matrix $S$ but use its full action. In this section we will outline their similarities and differences.

\begin{Remark} Since in the Lanczos method all residuals are linearly dependent, also the Lanczos method uses all information about the residual. However, Lanczos has no inner iteration. This will make sense from the point of view taken in Section \ref{sect3}. There we will show that for the choice $A_0=0$, also SPAM needs no inner iteration, and that this choice makes SPAM mathematically equivalent to Lanczos.
\end{Remark}

\subsubsection{The Riccati and the Jacobi-Davidson methods.}\label{sect2.5.1}

Given a Ritz pair $(\mu,u)$ with residual $\hat{r}\perp\VV$, generally each eigenvector of $A$ has a multiple that equals $u+t$, with $t\perp u$ a so-called \emph{orthogonal correction} to $u$. Indeed, let $X$ be such that $(u|X)$ is unitary. Then with $S=X^*AX$ and $\hat{r}=Xr$,
\be\label{c1} A(u|X) = (u|X) \left[\begin{array}{cl} \mu & r^* \\ r & S \end{array}\right],\ee
and the orthogonal correction $t$ equals $Xp$ where $p$ can be verified to satisfy the \emph{generalized algebraic Riccati equation}
\be \label{ric} (S-\mu I)p = -r + pr^*p.\ee
Transforming (\ref{ric}) back to the original basis, shows that $t$ solves
\be\label{ric1} t\perp u \und (I-uu^*)(A-\mu I)(I-uu^*)t = -\hat{r} + t\hat{r}^*t. \ee
In \cite{Bra}, solutions of (\ref{ric2}) were approximated by means of Rayleigh-Ritz projection in a $\ell$-dimensional subspace $\UU$ of $u^\perp$ and a suitable one was selected as expansion vector for $\VV$. This idea was intended as an enhancement of the Jacobi-Davidson method \cite{SlVo}. This method neglects the quadratic term $t\hat{r}^*t$ from (\ref{ric1}) and uses instead the unique solution $\hat{t}$ of
\be\label{ric2} \hat{t}\perp u \und (I-uu^*)(A-\mu I)(I-uu^*)\hat{t} = -\hat{r} \ee
to expand $\VV$. Jacobi-Davidson is simpler than the Riccati method in the sense that only a linear system for $\hat{t}$ needs to be solved. It can be interpreted as an {\em accelerated Newton method} \cite{SlVo2}. However, much more than the Riccati method, Jacobi-Davidson suffers from stagnation in case the term $t\hat{r}^*t$ from (\ref{ric1}) is not small. On the other hand, if $t\hat{r}^*t$ is small enough, Jacobi-Davidson converges quadratically, as one would expect from a Newton type method. This shows that one should be careful in proposing alternatives to the correction equation (\ref{ric2}). For instance, in \cite{GeSl}, the authors investigated the effect of solving the following alternative correction equation,
\be\label{ric4} \tilde{t}\perp V \und (I-VV^*)(A-\mu I)(I-VV^*)\tilde{t} = -\hat{r}. \ee
At first sight this seems to make sense, because it directly looks for a correction orthogonal to $\VV$. Also, the conditioning of the linear equation (\ref{ric4}) may be better than (\ref{ric2}) in case the search space contains good approximations of eigenvectors belonging to eigenvalues close to $\mu$. However, orthogonalizing $\hat{t}$ from (\ref{ric2}) to $\VV$ generally does not result in $\tilde{t}$ from (\ref{ric4}) and the price to pay is that $\|t-\tilde{t}\|$ is not of higher order, as is $\|t-\hat{t}\|$. Indeed, in \cite{GeSl} it is shown explicitly that, apart from some exceptional cases, the quadratic convergence of Jacobi-Davidson is lost, whereas in those exceptional cases, both expansions are equivalent. Numerical experiments in \cite{GeSl} confirm the above observations.

\begin{Remark} The method using the correction equation (\ref{ric4}) may at first sight also resemble SPAM. However, as we will see in the section to come, the use of $V$ in SPAM is of a different nature than in (\ref{ric4}), where, the correction is sought orthogonal to $\VV$. In SPAM, it is sought in the whole space and only orthogonalized to $\VV$ afterwards.
\end{Remark}

\subsubsection{Preconditioning in the Jacobi-Davidson method.}\label{osjd}

In the original paper \cite{SlVo} on the Jacobi-Davidson method, reprinted as \cite{SlVo3}, preconditioning is discussed as follows. Suppose that an approximation $A_0$ of $A$ is available. It is shown in \cite{SlVo} how to apply such a preconditioner to (\ref{ric2}), which is a linear equation, though not with system matrix $A$. To be explicit, since $(I-uu^*)\hat{t}=\hat{t}$, we have that
\be (A-\mu I)\hat{t} = -\vep u-\hat{r}. \ee
where $\vep$ is such that $\hat{t}\perp u$. Or equivalently, written as an augmented system,
\be\label{aug} \left[\begin{array}{cc}A-\mu I & u \\ u^* & 0\end{array}\right]\left[\begin{array}{c} \hat{t} \\ \vep \end{array}\right] = \left[\begin{array}{r} -\hat{r} \\ 0 \end{array}\right].\ee
Thus, an approximation $\hat{t}_0$ of $\hat{t}$, together with an approximation $\vep_0$ of $\vep$ can be obtained simply by replacing $A$ by $A_0$ in (\ref{aug}). The pair $\hat{t}_0,\vep_0$ can be computed as
\be\label{abc} \hat{t}_0 = -\vep_0 (A_0-\mu I)^{-1}u-(A_0-\mu I)^{-1}\hat{r}, \hdrie\mbox{\rm with } \vep_0 = -\frac{u^*(A_0-\mu I)^{-1}\hat{r}}{u^*(A_0-\mu I)^{-1}u}.\ee
This approximation is called a {\em one step approximation} in \cite{SlVo}. It was observed that setting $\vep_0=0$ in (\ref{abc}), the Davidson method \cite{Dav} results. With $A_0=A$, which corresponds to Jacobi-Davidson with full accuracy solution of the correction equation, (\ref{abc}) becomes
\be \hat{t} = -\vep(A-\mu I)^{-1} u - u, \ee
and since $\hat{t}$ is then orthogonalized to $u$, the method is mathematically equivalent to an accelerated shift and invert iteration that works with $(A-\mu I)^{-1}u$. It is argued, and demonstrated by experiments in \cite{SlVo}, that Jacobi-Davidson combines the best of those two methods. Of course, a natural next stage in preconditioning is to use the matrix
\be A_0^u =  \left[\begin{array}{cc}A_0-\mu I & u \\ u^* & 0\end{array}\right] \ee
as preconditioner {\em within the iterative method} that aims to solve (\ref{aug}). In each step of such a method one would need to solve a system with $A_0^u$. This can be done by solving two systems as in (\ref{aug})-(\ref{abc}) in the first step of the inner iteration. In each consecutive step, only one system of the form $(A_0-\mu I)z=y$ would need to be solved.

\subsubsection{One step approximation of the SPAM eigenproblem for $A_k$.}\label{SPAMJD}

In the SPAM method, the expansion vector for the Ritz Galerkin subspace in the outer iteration is a relevant eigenvector $v_k$ of $A_k$. In principle, any eigenvalue method can be used to compute an approximation for $v_k$, but observe that the starting point is as follows. In the outer iteration, we have just solved a $k\times k$ eigenproblem for $M=V^*AV$, and a Ritz pair $(\mu,u)$ with $\|u\|=1$ and with residual $\hat{r}=Au-\mu u$ has been selected. The matrix $A_k$ is now available, either explicitly or implicitly. Since $A_kV=AV$ and $V^*A_k=V^*A$, the Ritz pair $(\mu,u)$ for $A$ with respect to the current search space $\VV_k$ is also a Ritz pair for $A_k$. Thus we can exploit the fact that $\VV_k$ contains, by definition, good approximations of the relevant eigenvectors of $A_j$ with $j<k$, and use it as initial search space for a Ritz Galerkin method applied to $A_k$ to approximate $v_k$. Since $\VV_k$ is generally not a Krylov subspace, the Lanczos method is not a feasible candidate. The Jacobi-Davidson method is. The correction equation for the first step of the Jacobi-Davidson method in the inner iteration can be set up without any additional computations:
\be \label{jdspam} \left[\begin{array}{cc}A_k-\mu I & u \\ u^* & 0\end{array}\right]\left[\begin{array}{c} t_k \\ \vep_k \end{array}\right] = \left[\begin{array}{r} -\hat{r} \\ 0 \end{array}\right].\ee
Since quadratic convergence {\em in the outer iteration} cannot be expected even if an exact eigenvector of $A_k$ would be computed, we study the effect of applying only one iteration of Jacobi-Davidson in the inner iteration. This is also motivated by the fact that the initial search space $\VV_k$ for Jacobi-Davidson applied to $A_k$ may be relatively good and may result in quadratic convergence {\em in the inner iteration}.

\begin{Remark} If only one step of Jacobi-Davidson is applied, then after solving $t_k$ from (\ref{jdspam}), the new approximation $v$ for the eigenvector $v_k$ of $A_k$ would lie in the space $\VV_k\oplus\langle t_k\rangle$. It would not be necessary to actually compute this approximation, because
\be \VV_k\oplus\langle v\rangle \subset \VV_k\oplus\langle t_k\rangle.\ee
Thus, instead of computing the eigendata of a $(k+1)\times (k+1)$ matrix, the expansion of the Ritz-Galerkin space in the outer iteration can be done immediately with $t_k$.
\end{Remark}
SPAM, in which the eigenproblem for $A_k$ is approximated with one iteration of Jacobi-Davidson, we will refer to  as {\em one step SPAM}, abbreviated by SPAM(1). We will talk about {\em Full SPAM} if the eigenproblem for $A_k$ is solved to full precision.

\subsubsection{Comparing one step Jacobi-Davidson with SPAM(1).}\label{sect2.5.4}

SPAM(1) can best be compared with preconditioned Jacobi-Davidson with {\em one step approximation}, as described in Section \ref{osjd}. The only difference between the two method is that in iteration $k$ of SPAM(1) the preconditioner $A_k$ is used, whereas in one-step Jacobi-Davidson this is $A_0$. As such, SPAM(1) can be seen as an attempt to enhance this type of preconditioned Jacobi-Davidson. We will now investigate the effect of this attempt.

\begin{Lemma} Assume that $\VV_{k+1} = \VV_{k}\oplus \langle v\rangle$ with $v\perp\VV_{k}$ and $\|v\|=1$. Then
\be \label{expr} A_{k+1} = -vv^*Avv^* + Avv^* + vv^*A + (I-vv^*)A_{k}(I-vv^*). \ee
\end{Lemma}
{\em Proof. } By substitution of the defining relation for $A_{k}$. \hfill $\Box$

\begin{Corollary}\label{Cor9} Let $\VV_1$ be the span of the relevant eigenvector $u$ of $A_0$ and let $\mu=u^*Au$ and $\hat{r}=Au-\mu u$. Then the solution $t_1$ from the system
\be \left[\begin{array}{cc}A_1-\mu I & u \\ u^* & 0\end{array}\right]\left[\begin{array}{c} t_1 \\ \vep_1 \end{array}\right] = \left[\begin{array}{r} -\hat{r} \\ 0 \end{array}\right],\ee
in the first iteration of SPAM(1), to be used to expand $\VV_1$, coincides with the solution $t_0$ from the system
\be \left[\begin{array}{cc}A_0-\mu I & u \\ u^* & 0\end{array}\right]\left[\begin{array}{c} t_0 \\ \vep_0 \end{array}\right] = \left[\begin{array}{r} -\hat{r} \\ 0 \end{array}\right],\ee
solved in the first iteration of Jacobi-Davidson with one-step approximation using $A_0$.
\end{Corollary}
{\em Proof. } The linear system (\ref{jdspam}) for SPAM(1) with $k=1$ is equivalent to
\be t_1\perp u, \hdrie (I-uu^*)(A_1-\mu I)(I-uu^*)t_1 = -\hat{r},\ee
where $u$ is a unit vector spanning $\VV_1$. Substituting the expression (\ref{expr}) for $A_1$ immediately proves the statement. \hfill $\Box$

\smallskip
Thus, there is no difference in the first iteration. There is, however, a difference in further iterations. To study this difference we take a different viewpoint. Above, we approximated the SPAM inner eigenvalue problem by a linear correction equation which made it suitable for comparison with one step Jacobi-Davidson. The opposite viewpoint is also possible, which is to interpret the Jacobi-Davidson correction equation with one-step approximation as an {\em exact} correction equation of a {\em perturbed} eigenproblem.

\begin{Lemma} Given $u$ with $\|u\|=1$ and $\mu=u^*Au$ and $r=Au-\mu u$. Define for a given approximation $A_0$ of $A$ the matrix
\[  A_u :=-uu^*Auu^* + uu^*A + Auu^* + (I-uu^*)A_0(I-uu^*) \]
\be\label{sixteen} = u\mu u^* + ur^*+ru^*+(I-uu^*)A_0(I-uu^*). \ee
Then $(u,\mu)$ is a Ritz pair of $A_u$ in the one-dimensional span of $u$ with residual $r=A_uu-\mu u$ and with the equation
\be\label{sixb} (I-uu^*)(A_0-\mu I)(I-uu^*)t = -r\ee
as its exact (i.e., without preconditioning) Jacobi-Davidson correction equation.
\end{Lemma}
{\em Proof. }It is easily verified that $u^*A_uu=\mu$ and $A_u u-\mu u = r$ and thus $(\mu,u)$ is a Ritz pair for $A_u$ with residual $r$. Since moreover,
\be (I-uu^*)A_u(I-uu^*) = (I-uu^*)A_0(I-uu^*) \ee
its correction equation is precisely (\ref{sixb}), or equivalently, ({\ref{abc}). \hfill $\Box$

\bigskip
Note that $A_u$ from (\ref{sixteen}) is the subspace projected approximate matrix for the one dimensional subspace spanned by $u$. Now, with
\be \VV_k=\UU \oplus \langle u \rangle,\ee
where $\UU$ is the orthogonal complement of the span of the relevant Ritz vector $u$ in the current search space $\VV_k$, we have, similar as in (\ref{expr}) and (\ref{sixteen}), that
\be\label{seventeen} A_k = u\mu u^* + ur^* + ru^* + (I-uu^*)A_\UU (I-uu^*).\ee
Here, $A_\UU$ is the subspace projected approximated matrix corresponding to the subspace $\UU$. Now, the opposite viewpoint mentioned above, is the observation that in one step Jacobi-Davidson, the expansion vector is (an arguably good approximation of) the relevant eigenvector of $A_u$ in (\ref{sixteen}), whereas in SPAM, it is (an arguably good approximation of) the relevant eigenvector of $A_k$ in (\ref{seventeen}). Both matrices have now an appearance that is suitable for studying their differences and similarities.

\bigskip
The most important observation is that neither correction will lead to the unique correction that results in quadratic convergence in the outer iteration. Second, since both matrices $A_u$ and $A_k$ differ only in their restriction to the orthogonal complement $u^\perp$ of $u$, the difference in the methods they represent will be marginal if the residual is already small. Since, as already mentioned in Corollary \ref{Cor9}, not only at the start but also in the first iteration both methods coincide, the difference from the second iteration onwards will probably be very small, especially if $A_0$ provides a good initial approximation of the relevant eigenvector. Finally, even though $A_u$ uses less information of $A$ than $A_k$, it does use the {\em optimal} information in some sense. It may even be a disadvantage to use more information, because this involves approximations of eigenvectors orthogonal to the eigenvector of interest. The above observations will be tested, and confirmed by our numerical illustrations in Section \ref{sect4}.

\section{Selecting the matrix $A_0$ in the SPAM method}\label{sect3}
Here we study some of the effects of the choice for $A_0$ on the iterative approximation process in SPAM. We will assume for simplicity that the largest eigenvalue of $A$ is the target, although replacing $A$ by $-A$, we might as well have set the smallest eigenvalue of $A$ as a target.

\subsection{Some simple choices for $A_0$}\label{sect3.1}

It is instructive to study the consequences of the choice $A_0=0$. This generic parameter-free choice may seem dubious at first sight, but it is not. First note that since the start vector of SPAM is the relevant eigenvector of $A_0=0$, this is necessarily a random vector $v\in\CC^n$, with $\|v\|=1$. Thus, we set $\VV_1=\sspan\{v\}$. Write $\mu=v^*Av$ and $\hat{r} = Av-v\mu$. Then, with $V_\perp$ such that $(v|V_\perp)$ is orthogonal, and $r=V_\perp^*\hat{r}$ we find that
\be\label{A1} A = (v|V_\perp)\left[\begin{array}{cc}\mu & r^* \\ r & S\end{array}\right](v|V_\perp)^*,\ee
and consequently, replacing $S$ by the zero matrix, the next approximating matrix $A_1$ from (\ref{eq:ASPAM}) is defined by
\be \label{A2} A_1 = (v|V_\perp)\left[\begin{array}{cc} \mu & r^* \\ r & 0\end{array}\right](v|V_\perp)^*.\ee
As shown already in Proposition \ref{props1}, $A_1$ is a simple rank-two matrix, that on the basis defined by the columns of $(v|V_\perp)$ is of arrow-head type. It has two nontrivial eigenpairs $ A_1 w_\pm = \theta_\pm w_\pm$, where
\be\label{five} w_\pm  = \theta_\pm v+r \und \theta_\pm = \half \mu \pm \sqrt{\frac{1}{4}\mu^2+\|r\|^2}.\ee
Assuming that $A$ is positive definite will lead to the selection of $w_+$ for the expansion of the search space. Since $w_+$ is a linear combinations of $v$ and $r$, we find that
\be \VV_2 = \sspan\{v, r\} = \KK^2(A,v), \ee
and the two eigenvalue approximations computed in the outer loop of the SPAM method are the same as in the Lanczos method. This is, of course, not a coincidence.

\begin{Theorem}\label{prop4} If the goal is to find the largest eigenvalue of a positive definite matrix $A$, the SPAM method with $A_0=0$ is mathematically equivalent to the Lanczos method.
\end{Theorem}
{\em Proof. } Let $A_0=0$. Then the eigenvalues of $A_k$ in (\ref{eq:ASPAM}) are those of the $n\times n$ Hermitian arrowhead
\be \left[\begin{array}{cl} M & R^* \\ R & 0 \end{array}\right]. \ee
The Cauchy Interlace Theorem immediately gives that the largest $k$ of them are each larger than or equal to the corresponding eigenvalue of $M$. This assures that the eigenvector of $A_k$ that is selected for expansion is not from its null space but from its column span. From (\ref{stbasis}) we see that if $A_0=0$ this column span is the span of $V$ and $AV$. A simple induction argument shows that this span equals $\KK^{k+1}(A,v)$. \hfill $\Box$

\begin{Remark} If $A$ is indefinite and we would like to find the eigenvalue closest to zero, the choice $A_0=0$ would lead to expansion vectors from the null space of $A_0$, and the method would be worthless. As a solution, $A$ may be shifted by a suitable multiple $\alpha$ times the identity $I$ to make it positive semi-definite. Equivalently, instead of using $A_0=0$ we may choose $A_0=\alpha I$ as approximating matrix. In both cases, it is easy to verify that SPAM will still be equal to Lanczos.
\end{Remark}

The observant reader may have noticed that a peculiar situation has arisen. In the inner loop of SPAM, a better approximation of the largest eigenvalue of $A$ was computed than the Ritz values from the outer loop. In view of the philosophy of inner-outer iterations, this in itself is not out of the ordinary, but its computation did not require any additional matrix-vector multiplication with $A$, nor with an elaborate approximation $A_0$ of $A$. The following proposition, which uses the notation (\ref{eigs}) and (\ref{ritz}), makes this explicit.

\begin{Proposition}\label{prop1} With $\theta_+$ as in (\ref{five}) we have that $\mu \leq \|Av\| \leq \theta_+$. Moreover, if $A$ is positive semi-definite, we find that $\theta_+ \leq \lambda_+$, where $\lambda_+$ is the largest eigenvalue of $A$.
\end{Proposition}
{\em Proof. } Because  $r\perp v$ and $Av-\mu v = r$, by Pythagoras' Theorem we have
\be\label{three} \mu^2\|v\|^2 + \|r\|^2 =  \mu^2 + \|r\|^2 = \|Av\|^2, \ee
hence $\mu\leq \|Av\|$. Squaring $\theta_+$ from (\ref{five}) gives
\be \theta_+^2 = \|r\|^2 + \half \mu^2 + \mu \sqrt{\frac{1}{4}\mu^2+\|r\|^2},\ee
which together with (\ref{three}) shows that $\|Av\| \leq \theta_+$. Since $S$ in (\ref{A1}) is positive definite whenever $A$ is, combining (\ref{A1}) and (\ref{A2}) with Weyl's bound yields $\theta_+ \leq \lambda_+$.\hfill $\Box$\\[2mm]
The key issue here is that the inner eigenvalue approximations are much related to the so-called \emph{harmonic Ritz values} \cite{Bea,SlEs} of $A$. Indeed, assuming that $A$ itself is positive semi-definite, these are the $k$ positive eigenvalues of the at most rank-$2k$ matrix
\be\label{harm} \tilde{A}_k = \left[\begin{array}{cl} M & R^* \\ R & T\end{array}\right], \hdrie\mbox{\rm where } \hdrie T= RM^{-1}R^*. \ee
They can be computed without additional matrix vector multiplications with $A$. Note that harmonic Ritz values are usually introduced as the reciprocals of the Rayleigh-Ritz approximations of $A^{-1}$ in the space $A\VV$. It is well known that for positive semi-definite matrices $A$, the harmonic Ritz values are better approximations of the larger eigenvalues of $A$ than the standard Ritz values. We provide the short argument in Lemma \ref{th1}. See also \cite{Bea}.

\begin{Proposition} The matrix $\tilde{A}_k$ can be decomposed as
\be \left[\begin{array}{cl} M & R^* \\ R & T \end{array}\right] = \left[\begin{array}{c} M \\ R\end{array}\right] M^{-1} \left[\begin{array}{cc} M & R^*\end{array}\right]. \ee
The blocks
\be \left[\begin{array}{c} M \\ R \end{array}\right] \und \left[\begin{array}{c} -M^{-1}R^* \\ I \end{array}\right]\ee
span the range and null space, respectively, and the nonzero eigenvalues are the eigenvalues of the $k\times k$ matrix
\be UM^{-1}U^* \hdrie \mbox{\rm where }\hdrie  \left[\begin{array}{c} M \\ R \end{array}\right] = QU \ee
is a QR-decomposition. In particular, those eigenvalues are positive.
\end{Proposition}
{\em Proof. } The statements are all easy to verify. The positivity of the $k$ nonzero eigenvalues follows from Sylvester's Theorem of Inertia. \hfill $\Box$\\[2mm]
The $k$ eigenpairs $(\theta_j,w_j)$ of $\tilde{A}_k$ in (\ref{harm}) with positive eigenvalues we label as
\be 0 < \tilde{\theta}_{k} \leq \tilde{\theta}_{k-1} \leq \dots \leq \tilde{\theta}_{2} \leq \tilde{\theta}_1. \ee
The proposition shows that they can be computed by solving a $k\times k$ eigenproblem and that no additional matrix-vector products with $A$ are needed. We can now easily prove the following bounds. See also \cite{Bea,Par}.

\begin{Lemma}\label{th1} For all $j\in\{1,\dots,k\}$, we have that
\be\label{seven} \mu_j \leq  \tilde{\theta}_j \leq \lambda_j. \ee
\end{Lemma}
{\em Proof. } The left inequalities follow from the Cauchy Interlace Theorem applied to $\tilde{A}_k$. Now, with $\hat{A}$ as in (\ref{eq:Asubsp}), recognizing the Schur complement $\hat{A}/M$ shows that,
\be  \hat{A} - \tilde{A}_k = \left[\begin{array}{cl} 0 & 0^* \\ 0 & \hat{A}/M \end{array}\right]\ee
is positive semi-definite, hence the right inequalities follow from Weyl's bound.\hfill $\Box$

\bigskip
Observe that, as was the case with the choice $A_0=0$, assuming that from the start of the SPAM method the eigenvector $w_1$ belonging to $\theta_1$ was selected for expansion, the eigenvectors $w_j$, called the \emph{harmonic Ritz vectors}, lie in the column span of $AV$ and hence, in $\KK^{k+1}(A,v)$.

\bigskip
Thus, even though we may have improved the inner loop eigenvalue approximations, the SPAM method is still equal to the Lanczos method. It does give, however, valuable insight into SPAM: Lanczos results from the choice $A_0=0$, even after modification of $A_k$ into the positive semi-definite matrix $\tilde{A}_k$. In order to get a method different than Lanczos, we should use less trivial approximations that are based on the structure and properties of $A$ itself, while aiming to retain similar inequalities as in Lemma \ref{th1}. For this, we would need the matrices $A-A_k$ to be positive semi-definite. We will investigate this in the following sections.

\subsection{One-sided approximations}

Having seen that the trivial choice $A_0=0$, even after after a correction that turns all approximate matrices positive semi-definite, will generally lead to the Lanczos method, we now turn our attention to {\em approximations from below}, by which we mean $A_0$ such that $A-A_0$ is positive semi-definite.

\begin{Lemma}\label{lembelow} If $A_0$ approximates $A$ from below, then so does each matrix $A_k$.
\end{Lemma}
{\em Proof.} Combining (\ref{eq:Asubsp}) and (\ref{acht}) with (\ref{eq:ASPAM}) we see that for all $x\in\mathbb{C}^{n}$,
\begin{equation}\label{eq33}
x^{*}(\hat{A}-\hat{A}_k)x=
x^{*}\left[\begin{array}{cc}
    0 & 0^{*} \\
    0 & S-\tilde{S}
\end{array}
\right] x= x^{*}V_\perp^*(A-A_0)V_\perp^{\phantom{*}} x\geq 0
\end{equation}
where the last inequality holds because $A-A_0$ is positive semi-definite. And thus, $\hat{A}-\hat{A_k}$ is positive semi-definite, and hence, so is $A-A_k$. \hfill $\Box$

\bigskip
By Proposition \ref{props1}, $A_k-A_{k-1}$ is an indefinite rank-$2$ matrix, hence it will generally not be true that $A_{k-1}$ approximates $A_k$ from below.

\begin{Lemma}\label{lem:psd} The following inequalities are valid generally,
\be \theta_{j+1} \leq \mu_j \leq \theta_j \hdrie \mbox{\rm for all } j\in\{1,\dots,k\},\ee
whereas, if $A_0$ approximates $A$ from below, additionally
\be \theta_j \leq \lambda _j \hdrie \mbox{\rm for all } j\in\{1,\dots,n\}.\ee
\end{Lemma}
{\em Proof. } The first set of inequalities applies because $A_k$ has the same Ritz values as $A$. See also Section \ref{sect2.4.1}. It is well known that the Ritz values interlace the exact eigenvalues \cite{Par}. Since $A-A_k$ is positive semi-definite due to Lemma \ref{lembelow}, the equality $A_k+(A-A_k)=A$ together with Weyl's bound \cite{Par} proves the second set of inequalities.\hfill $\Box$.

\bigskip
Lemma \ref{lem:psd} shows that if $A_0$ approximates $A$ from below, the approximations for the larger eigenvalues of $A$ that are produced in the inner iteration, will never be worse than the ones obtained in the outer iteration. Moreover, they will never be larger than the corresponding exact eigenvalues. Thus, it indeed makes sense to expand the search space with the eigenvector that is computed in the inner iteration. Question that remains is how to obtain matrices $A_0$ that approximate $A$ from below.

\subsubsection{Algebraic construction of approximations from below.}\label{algebraic}

Clearly, for any positive definite matrix $H$ we have that $A_0=A-H$ approximates $A$ from below, even though $A_0$ itself may not be positive definite. The problem is of course how to choose $H$ such that $A_0$ is close to $A$ in an appropriate sense, while its action is considerably less expensive than that of $A$.

\bigskip
If $A$ itself is a positive definite matrix, a purely algebraic option is at hand. Given an index set $\mathcal{I}\subset \{1,\dots,n\}$ of cardinality $m$, let $E_\II$ be the matrix with the standard basis vectors $e_i,i\in \II$ as columns. Then set
\be A_0 = A - H, \hdrie \mbox{\rm where }\hdrie H = E_\II^{\phantom{*}} E_\II^*AE_\II^{\phantom{*}} E_\II^{*}.  \ee
The matrix $H$ is the result of a Rayleigh-Ritz procedure with as search space the column span of $E_\II$. For a randomly chosen index set, this space has no a priori relation with $A$ and thus $H$ is probably a relatively poor approximation of $A$ in comparison with, for example, a Krylov subspace approximation.

\begin{Remark} In this particular situation it is an advantage if $H$ does not approximate $A$ very well, because $A_0$ should be close to $A$, not $H$. Notice also that $A_0$ has zero entries at positions $(i,j)$ for all $i,j\in\II$, and is thus always more sparse than $A$. A priori knowledge of $A$ may lead to a more sophisticated choice of the index set(s) $\II$. If the goal is to approximate the largest eigenvalues of $A$, the index set $\II$ could be chosen such that the smallest diagonal entries of $A$ are selected to put in $H$. Consequently, $A_0$ will share with $A$ the largest diagonal entries, and this may increase its approximation quality. This is illustrated in Figure \ref{figure2}.
\end{Remark}

\begin{figure}
\centering
\includegraphics[width=6.5cm, height=7cm]{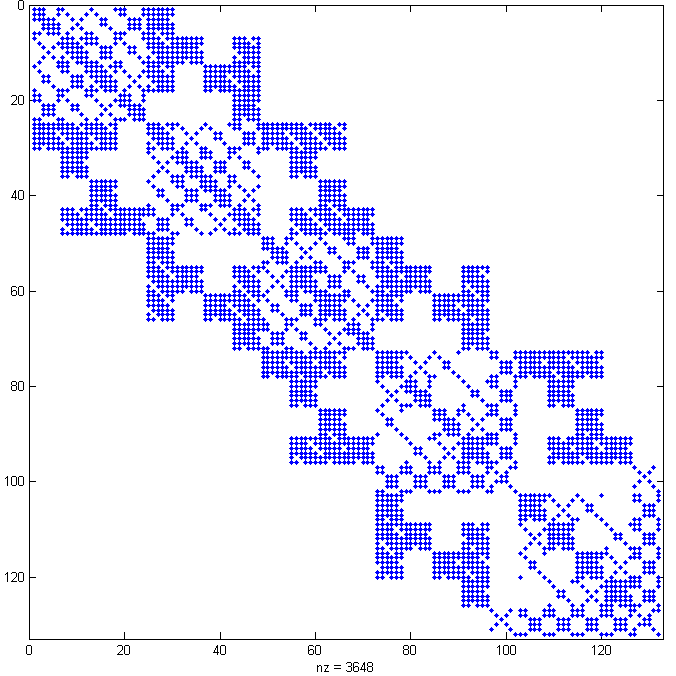} \includegraphics[width=6.5cm, height=7cm]{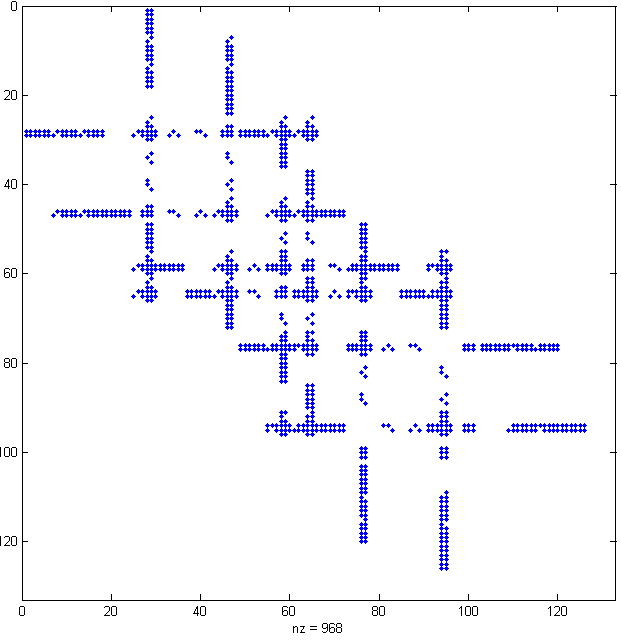}
\caption{Approximating a matrix $A$ from structural engineering from below by a sparser matrix $A_0$ by subtracting from $A$ a definite matrix $H$; the sparsity plot of $A$ (with $3648$ nonzero entries) is on the left and of $A_0$ (with $968$ nonzero entries) on the right. See the experiments with this pair of matrices in Section \ref{sect4.4}.} \label{figure2}
\end{figure}

\begin{Remark}{\rm Notice that $\rank(A_0)\leq 2(n-m)$. Especially for large $m$ this may greatly simplify the computation of the eigendata of $A_0$ and of $A_k$ for small values of $k$.}
\end{Remark}

\begin{Remark} Although $A$ itself is positive definite, $A_0$ generally is not. It can be made positive definite by adding a Schur complement at the position of the zero block, as was done in (\ref{harm}). Since the Schur complement block is the smallest correction to $A_0$ that makes it positive semi-definite, the result would still approximate $A$ from below. However, its computation involves the evaluation of the inverse of an $(n-m)\times(n-m)$ matrix. It is not clear if the additional computational effort is well spent. Also without the Schur complement, $A_0$ has as many positive as negative eigenvalues and is likely to perform better than $A_0=0$.
\end{Remark}

The rather crude approach that we just described is easy to implement, and it can be applied in the context of the {\em multilevel version} of SPAM mentioned in Section \ref{sect2.4.2}: the solution of the inner iteration is done using SPAM itself, with an approximation of $A_0$ that is based on a larger index set than the one that was used to construct $A_0$ itself.

\subsubsection{Natural construction of approximations from below.}\label{natural}

A situation in which approximations from below are naturally available is the setting of discretized partial differential equations including an elliptic term, either by the finite difference method or the finite element method \cite{CiLi}. Then removing the positive definite discrete elliptic operator, either completely or partially, from the total discrete operator, results in an approximation from below. Indeed, let $\Omega\subset\RR^3$ be a domain, and consider as an example the problem of finding eigenmodes for the linear operator $L$, defined by
\be\label{pdv} L(u) = \lambda u \hdrie\mbox{\rm where }\hdrie L(u) = -\varepsilon {\rm div}(\mathcal{K} \nabla u) + cu, \und u=0 \hdrie\mbox{\rm on } \partial\Omega.\ee
Here $\varepsilon>0$ is a parameter, $\mathcal{K}:\Omega\rightarrow \mathbb{M}^{3\times 3}(\RR)$ maps into the symmetric positive definite matrices, and $c\in C(\overline{\Omega})$ is nonnegative. Discretizing this equation with the finite difference method will lead to an algebraic eigenvalue problem of the form $\varepsilon K x + M x = \xi x.$ The matrix $K$ that represents the discretized diffusion is positive definite. Although sparse, it will generally have more fill-in than the diagonal matrix $M$ that represents the reaction term, and, if $\varepsilon$ is small enough, have smaller entries. Thus, the total discretized operator $A=K+M$ has $A_0=M$ as a candidate for the approximation matrix: its action is cheaper than the action of $A$ and $A-A_0=K$ is positive definite. A similar strategy can also be employed in the finite element method, when so-called \emph{mass lumping} is used in the assembly of the mass matrices.

\begin{Remark} In this context, the algebraic method can be used to approximate the smallest eigenvalue of $A$ by applying it to $\alpha I-A$ with $\alpha$ such that $\alpha I-A$ is positive semi-definite. Even though the largest eigenvalue usually has no physical relevance, together they would provide good estimates for the condition number of the $A$, which is indeed of interest.
\end{Remark}

\subsection{Cutting off the bandwidth} \label{sect3.5}
Here we describe an obvious choice of approximating matrices that was used in \cite{ShWaTiMi} to illustrate the effectiveness of their method. It concerns symmetric banded matrices. Apart from the approach that we will now describe, in the numerical illustrations of Section 4 we also intend to apply our algebraic approximations from below to these matrices.\\[2mm]
Given $0\leq \epsilon\leq 1$ and $1\leq q\leq n-1$, define a matrix $A=(A_{ij})$ by
\begin{equation}\label{eq:BandMat}
    A_{ij}=\begin{cases}
            i & \text{if $i=j$; } \\
            \epsilon^{q} & \text{if $1\leq |i-j|\leq q$;}\\
            0 & \text{otherwise.}
     \end{cases}
\end{equation}
In \cite{ShWaTiMi} it is proposed to choose as approximating matrix for $A$ a matrix $A_0$ of the same type with a smaller half-bandwidth $\tilde{q}<q$. For instance,
\begin{equation}\nonumber
    A=\left[\begin{array}{cccc}
      1 & \epsilon & \epsilon^{2} &  \\
     \epsilon & 2 & \epsilon &\epsilon^{2} \\
     \epsilon^{2} & \epsilon & 3 & \epsilon \\
      & \epsilon^{2} & \epsilon & 4
\end{array}\right],\hspace{5mm}\mbox{\rm then with }
    A_0=\left[\begin{array}{cccc}
      1 & \epsilon &  &  \\
      \epsilon & 2 & \epsilon & \\
       & \epsilon & 3 & \epsilon \\
       &  & \epsilon & 4
\end{array}\right].
\end{equation}
For each eigenvalue $\theta$ of $A_0$ there is an eigenvalue $\lambda$ of $A$ with
\begin{equation}\label{bound}
        |\lambda-\theta|\leq \frac{\epsilon^{q_0+1}-\epsilon^{q+1}}{1-\epsilon}\sqrt{n}.
    \end{equation}
Indeed, the difference $A-A_0$ is zero except for the bands $q_0+1$ to $q$. Each non-zero row contains at most the numbers $\epsilon^{q_0+1}$ to $\epsilon^{q}$, and (\ref{bound}) follows from the the Bauer-Fike \cite{BaFi} theorem and the finite geometric series sum formula.
Thus, for values of $\epsilon$ small enough, $A_0$ may be a good candidate to approximate $A$ even though it is generally not an approximation from below. Nevertheless, its eigenvalues are close to the eigenvalues of $A$ and its action is cheaper than the action of $A$. The number of floating point operations required to compute $Av$ is approximately twice as much as for an approximation $A_0$ having half the bandwidth of $A$. In other words, each two matrix-vector products $A_0$ are approximately equally costly as a single product with $A$.

\bigskip
Cutting off the bandwidth in this fashion makes sense especially if the decay of the size of the off-diagonal entries in relation to their distance to the main diagonal is quick enough. Apart from the example above, this is the case in many applications where the boundary element method \cite{WrAl} is used to approximate integral equations. Notice that for such applications also the approach from Section \ref{algebraic} can be applied, resulting in both sparse and low rank approximating matrices $A_0$.

\section{Numerical illustrations}\label{sect4}

In the previous sections, we have compared SPAM with both the Jacobi-Davidson method and the Lanczos method. We have also studied ways to define an appropriate approximate matrix $A_0$ in the context of the SPAM method. Since the choice $A_0=0$ will effectively lead to the Lanczos method, our starting point in the upcoming numerical illustrations will be to consider SPAM as a boosted version of the Lanczos method. This is, of course, particularly justified if $A_0$ is an approximation of $A$ from below. As discussed in Section \ref{Sect2.3}, SPAM can also be considered as an attempt to enhance the Jacobi-Davidson method with one-step approximation as preconditioning. Therefore, we will also present a comparison of SPAM and this version of Jacobi-Davidson. We end this section with a discussion of the numerical results.

\subsection{Objectives}
First, we list the methods and abbreviations that we use to describe them.
\begin{itemize}
\item Lanczos (see Section \ref{sect-Lanczos});
\item JD($\ell$): Jacobi-Davidson using $\ell$ steps of MinRES to approximate the solution of the exact correction equation (\ref{aug}) in augmented form;
\item JD(1,$\ell$): Jacobi-Davidson with one step approximation as preconditioner (see Section \ref{osjd}) using the matrix $A_0$, and $\ell$ steps of MinRES to approximate the solution of the correction equation (\ref{aug}), with $A$ replaced by $A_0$;
\item Full SPAM: each eigenproblem for $A_k$ solved in full precision;
\item SPAM(1): eigenproblem for $A_k$ approximated with one step of the Jacobi-Davidson method, correction equation (\ref{jdspam}) in augmented form solved to full precision (see Section \ref{SPAMJD});
\item SPAM(1,$\ell$): using $\ell$ steps of MinRES to approximate the solution of the correction equation for SPAM(1).
\end{itemize}
\begin{Remark} To minimize the size of the legends in the pictures, we sometimes write LZS for Lanczos, FSP for Full SPAM, SP(1) for SPAM(1), S12 for SPAM(1,2), JD13 for JD(1,3), etcetera.
\end{Remark}
In the experiments, we will give illustrations of the following aspects of SPAM.
\begin{itemize}
\item When a nonzero approximation $A_0$ of $A$ from below is used, less outer iterations of Full SPAM are needed to arrive close to the dominant eigenvalue of $A$ than with the choice $A_0=0$, which is equivalent to the Lanczos method with a random start vector.
\item Even if the Lanczos method is started with the same approximation of the dominant eigenvector of $A_0$ as Full SPAM, Full SPAM method will still outperform Lanczos in terms of the number of outer iterations.
\item Also for other eigenvalues, Full SPAM outperforms Lanczos; this may be expected because in Full SPAM the Ritz Galerkin subspace will be forced in the direction of the appropriate eigenvector in every iteration. In Lanczos, this is done only by means of the start vector. Of course, Lanczos allows efficient implicit restart strategies \cite{Sor}, but also Full SPAM may be restarted. We feel that the comparison would become diffuse and thus we refrain from incorporating restarts.
\item We investigate the effect of approximating the desired eigenvector of $A_k$ with just one step of the Jacobi-Davidson method, i.e., we will be comparing Full SPAM with SPAM(1).
\item SPAM(1,$\ell$) will be compared with JD(1,$\ell$), both started with the same initial vector, i.e., the relevant eigenvector of $A_0$; i.e., both methods will spend the same number $\ell$ of matrix-vector products in their inner iteration, where in SPAM(1,$\ell$), the matrix will be $A_k$, and in JD(1,$\ell$), the matrix will be $A_0$. From the viewpoint of this paper, this is the comparison that is of most interest. Not only is JD(1,$\ell$) the closest related to SPAM(1,$\ell$), the difference between the two is solely the fact that the action of $A$ from the outer loop is taken into the inner iteration of SPAM(1,$\ell$), whereas in JD(1,$\ell$), this is not done. See the discussion in, and particularly at the end of Section \ref{sect2.5.4}.
\item Finally, we compare SPAM(1,$\ell$) with JD($\ell$). This is perhaps the comparison that the authors of SPAM \cite{ShWaTiMi} had in mind: in the inner iteration of JD($\ell$) the original matrix $A$ is used, whereas in SPAM(1,$\ell$) it will be $A_k$.
\end{itemize}

We will comment on the computational costs of an inner iteration in comparison to having no such costs (as in Lanczos) or the full costs (Jacobi-Davidson), although these costs may depend very much on the specific problem and the available approximations and even on hardware parameters like available memory.

\subsection{Lanczos versus full SPAM: Reaction-Diffusion problem, various eigenvalues} \label{sect4.1}
In this section we will compare Lanczos with Full SPAM. Our comparison is, for the time being, only in terms of the number of outer iterations. A first naive comparison uses a random start vector for Lanczos, but from then onwards, we will start Lanczos with the appropriate eigenvector of $A_0$. The approximate matrix $A_0$ will be contructed using both approaches described in Sections \ref{algebraic} and \ref{natural}. For this, we discretized the one-dimensional version of (\ref{pdv}) on $\Omega=[0,1]$ using finite differences with grid size $h=1/33$ with the parameters $\varepsilon = \frac{1}{33^2}, \mathcal{K} = 1$, and $c(x) = x(1-x){\rm e}^{3x}$. The resulting $32\times 32$ algebraic eigenproblem is of the form $Ax=\lambda x$, where $A=D+R$ is the sum of the tridiagonal discretized diffusion $D$ and the diagonal discretized reaction $R$. With the approach from Section \ref{natural} we approximate the largest eigenvalue and with the approach from Section \ref{algebraic} the smallest eigenvalue.

\paragraph{Natural approximation from below} The left picture of Figure \ref{Figure4} illustrates the typical convergence of the Lanczos method with a random start vector. We display the $k$ Ritz values at outer iteration $k$ as circles above the value $k$ of the iteration number on the horizontal axis.
Due to the interlacing property, eigenvalue approximations ``converge'' to either side of the spectrum, and this is emphasized by the connecting lines between Ritz values for different values of $k$, both upwards and downwards. In view of Theorem \ref{prop4} we could also say that this picture belongs to SPAM with choice $A_0=0$. In the right picture of Figure \ref{Figure4}, we show in a similar fashion the convergence of SPAM, using $A_0=R$ as the approximate matrix, as suggested in Section \ref{natural}. We see that the convergence towards the largest eigenvalues is stimulated. The costs for this faster convergence is solving a diagonal plus rank $2k$ eigenproblem in iteration step $k$. There exist efficient methods for such eigenproblems based on the secular equation and Newton's method.

\paragraph{Algebraic approximation from below} We also tested the algebraic approach from Section \ref{algebraic} to construct approximations from below. We created a rank-$12$ approximation $A_0$ of $A$ based on the largest diagonal elements of $A$. Full SPAM and Lanczos were started with the dominant eigenvector of $A_0$ in order to approximate its dominant eigenvalue. The leftmost picture in Figure \ref{Figure5} shows that the incorporation of an approximation $A_0$ into Full SPAM has an effect that carries beyond that of creating a better start vector for Lanczos. In the rightmost picture of Figure \ref{Figure5}, the same was done for the matrix $6I-A$, which is positive definite and has the smallest eigenvalue of $A$ as dominant eigenvalue. Also here, SPAM outperforms Lanczos in terms of the number of outer iterations. In the middle two pictures of Figure \ref{Figure5}, we plotted the absolute error in the second and fifth largest eigenvalue of $A$. The results are shown from iteration $2$ and $5$ onwards, respectively, because the Ritz-Galerkin method produces approximations of the second and fifth largest eigenvalues from that iteration onwards. Again, Full SPAM clearly outperforms Lanczos in both cases. Note that, when using Full SPAM to approximate the $p$-th largest eigenvalue, the Ritz-Galerkin subspace in the outer iteration is in each step expanded with the eigenvector belonging to the $p$-th largest eigenvalue of $A_k$. For a fair comparison, we also started Lanczos with the eigenvector of $A_0$ belonging to the $p$-th largest eigenvalue.

\begin{figure}[ht]
\includegraphics[width=7cm]{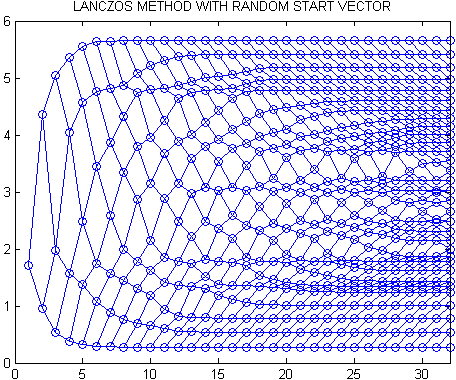} \includegraphics[width=7cm]{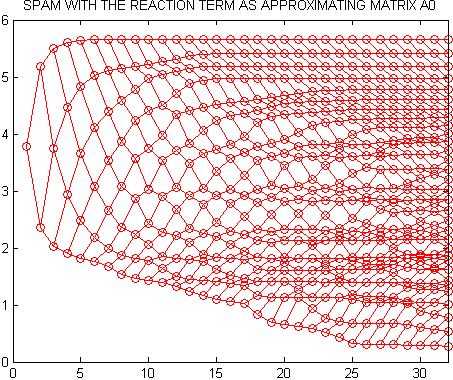}
\caption{Lanczos method (left) with a random start vector, versus SPAM (right) with the discretized reaction term as approximating matrix $A_0$, and the largest eigenvalue of $A$ as target.}\label{Figure4}
\end{figure}
\begin{figure}[ht]
\includegraphics[width=3.5cm, height=6cm]{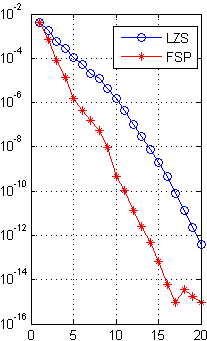} \includegraphics[width=3.5cm, height=6cm]{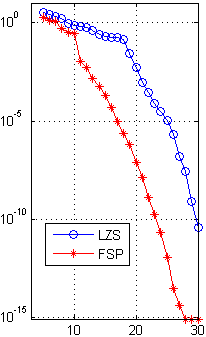}
\includegraphics[width=3.5cm, height=6cm]{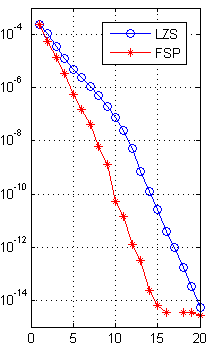} \includegraphics[width=3.5cm, height=6cm]{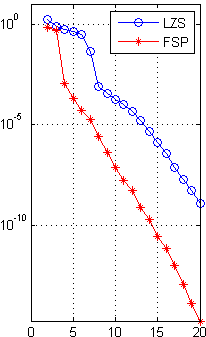}
\caption{From left to right: Lanczos and Full SPAM approximating the largest, the second, the fifth, and the smallest eigenvalue of $A$, using algebraic rank-$12$ approximation from below (for the smallest eigenvalue, we applied both the methods to $6I-A$). Lanczos and Full SPAM  used in each experiment the same start vector.}\label{Figure5}
\end{figure}

\subsection{Lanczos versus Full SPAM and SPAM(1): Banded matrices} \label{sect4.3}
In this section we not only compare Lanczos with Full SPAM, but also with SPAM(1), by which we mean that the eigenproblem for $A_k$ in Full SPAM is approximated using one step of the Jacobi-Davidson method, as explained in Section \ref{SPAMJD}. The correction equation that results is still solved to full precision. For our of experiments, we took the banded matrix from Section \ref{sect3.5} of size $32\times 32$ and with $q=5$ and $\epsilon=0.5$. In \cite{ShWaTiMi}, it was suggested to take for $A_0$ the matrix $A$ with a number of its outer diagonals put to zero. We repeated that experiment with $A_0$ equal to the diagonal of $A$, such that here too, $A_k$ is diagonal plus a rank-$2k$ perturbation. The comparison between Lanczos, Full SPAM and SPAM(1) is depicted in the left graph in Figure \ref{Figure7}. This comparison is initially in favor of Lanczos. This may be due to the fact that the difference $A-A_0$ is indefinite: it has only eleven nonnegative eigenvalues. In the middle left picture we took the tridiagonal part as approximation. This approximation is indeed positive definite and gives better results. Taking the simple algebraic approximation $A_0$ from below from Section \ref{algebraic} based on the three largest diagonal entries of $A$, which is of rank $6$, gives comparable results, but its low rank and higher sparsity make this choice more interesting. In the right graph in Figure \ref{Figure7}, the largest eigenvalue of $\alpha I-A$ was approximated with the approximation from below that kept the largest three diagonal entries of $\alpha I-A$, and thus the smallest three diagonal entries of $A$. In all cases, the positive effect of incorporating an approximation $A_0$ goes beyond delivering a good start vector for Lanczos. Also in all cases, there is virtually no difference between Full SPAM and SPAM(1).

\begin{figure}[ht]
\centering
\includegraphics[width=3.5cm, height=6cm]{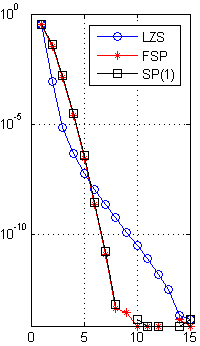}
\includegraphics[width=3.5cm, height=6cm]{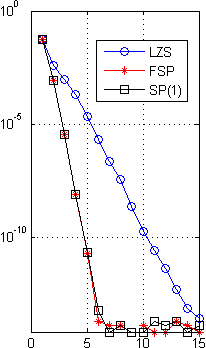}
\includegraphics[width=3.5cm, height=6cm]{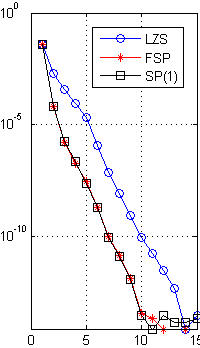}
\includegraphics[width=3.5cm, height=6cm]{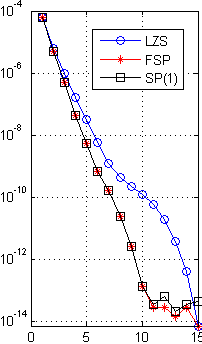}
\caption{Lanczos versus Full SPAM and SPAM(1) with diagonal (left), tridiagonal (middle left), and with algebraic rank-$6$ approximation from below (both middle right and right). In the three leftmost pictures, the target was the largest eigenvalue, at the right it was the smallest eigenvalue (i.e., the largest of $\alpha I-A$).}\label{Figure7}
\end{figure}

\subsection{Lanczos versus Full SPAM and SPAM(1): matrices from structural engineering} \label{sect4.4}
As a next set of experiments, we took some matrices from the Harwell-Boeing collection of test matrices. They have their origin in the area of structural engineering and are called {\em bcsstk04, bcsstk07} and {\em bcsstk10} and have respective sizes $132\times 132, 420\times 420$ and $1024\times 1024$. As approximating matrix $A_0$ we took approximations from below keeping respectively the largest $12, 20$ and $180$ diagonal entries. Recall that in Figure \ref{figure2} we displayed the sparsity plots of $A$ and $A_0$ for bcsstk04.mtx. As was the case for the banded matrix in the previous section, Full SPAM and SPAM(1) behave virtually the same, and need less outer iterations than the Lanczos method to arrive at a given accuracy.

\begin{figure}[ht]
\centering
\includegraphics[width=4.7cm, height=6.5cm]{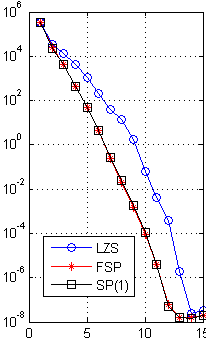} \includegraphics[width=4.7cm, height=6.6cm]{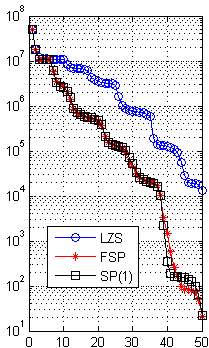}
\includegraphics[width=4.7cm, height=6.4cm]{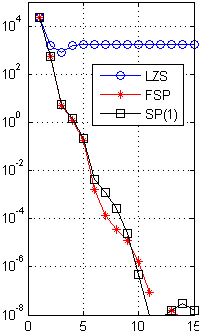}
\caption{Lanczos versus Full SPAM and SPAM(1) for bcsstk04, bcsstk07 and  bcsstk10, with approximating matrices from below. All three methods used the same start vector.}\label{Figure8}
\end{figure}

\paragraph{Conclusions} The main goal of the experiments so far was, first of all, to investigate if {\em Full} SPAM is competitive with the Lanczos method in terms of the number of iterations in the outer loop. Solving the eigenproblems for $A_k$ in the inner loop, even if this would not be done to full precision, makes SPAM always more expensive than Lanczos. The experiments in Figures \ref{Figure4}, \ref{Figure5}, \ref{Figure7} and \ref{Figure8} show that this is the case for different matrices $A$ and different types of approximations $A_0$. Not only for the largest, but also for other eigenvalues of $A$. Secondly, we have illustrated that the use of one step of the Jacobi-Davidson method to approximate the eigenvalue problem for $A_k$ in the inner iteration for SPAM hardly influences its behavior. We will now investigate what happens if the linear system in SPAM(1) is approximated with a small number of steps of the Minimal Residual method, MinRES \cite{PaSa}.

\subsection{Lanczos versus SPAM(1,$\ell$): approximating the correction equation of SPAM(1) using MinRES}

In this section we investigate the use of $\ell$ iterations of the MinRES method \cite{PaSa} for approximating the solution of the Jacobi-Davidson correction equation (\ref{jdspam}) for the eigenproblem for $A_k$ in the inner iteration of SPAM(1). Each iteration of MinRES requires one matrix vector product with $A_k$, which, for small values of $k$, will be approximately as costly as a matrix vector product with $A_0$. The initial vector of all methods to which the comparison is applied, i.e.,  the eigenvector of interest of $A_0$, will still be computed in full precision. The resulting method is abbreviated by SPAM(1,$\ell$).

\begin{Remark} In SPAM(1,1), MinRES produces an approximation of the Jacobi-Davidson correction equation (\ref{jdspam}) for $A_k$ in the one-dimensional Krylov subspace with the right-hand side of that correction equation as start vector. Since this is the current residual from the outer iteration, the expansion is the same as for the Lanczos method. We will therefore not display the convergence curves of SPAM(1,1).
\end{Remark}
In the light of the previous remark, it is reasonable to expect that SPAM(1,$\ell$) will represent a transition between Lanczos and SPAM(1). Thus, for reference, in the experiments to come, we displayed the convergence graphs of Lanczos and SPAM(1) in solid black lines without any additional symbols. The experiments concern four of the situations that we have already studied. First, we approximated the smallest eigenvalue of the reaction-diffusion problem. For the other three experiments we approximated the largest eigenvalue: in the second, we took the banded matrix with low rank approximation from below, and in the third and fourth the matrices bcsstk07 and bcsstk10 with the respective approximations from the previous section. The results are displayed in Figure \ref{Figure9} and confirm the expectations. Even for $\ell=2$ and $\ell=3$, SPAM(1,$\ell$) resembles SPAM(1) much more than it resembles Lanczos. It depends, however, very much on the actual application if the gain in the number of iterations is not undone by the costs of $\ell$ steps of MinRES per outer iteration with the matrix $A_k$. For instance, in the banded matrix example (second picture in Figure \ref{Figure9}), the matrix $A$ itself has $322$ nonzero entries and is of full rank, whereas $A_0$ only has $33$ nonzero elements and is of rank $6$, and especially when $k$ is small, the action of $A_k$ will not be very much more expensive than the action of $A_0$. This, however, brings us to the next question that we would like to investigate, which is, if using $A_k$ in the $k$-th inner iteration instead of $A_0$ all along, is going to make any difference, because this is what distinguishes SPAM from Jacobi-Davidson with one-step approximation as preconditioner. As argued in Section \ref{sect2.5.4}, if SPAM is going to to better, then probably not by very much.

\begin{figure}[ht]
\includegraphics[width=3.5cm, height=6cm]{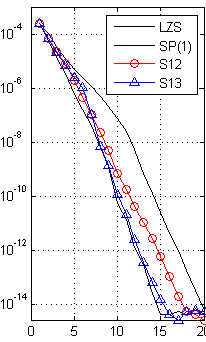}
\includegraphics[width=3.5cm, height=6cm]{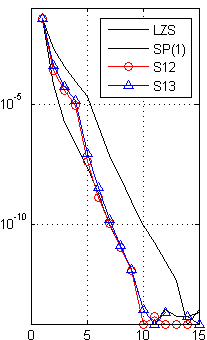}
\includegraphics[width=3.5cm, height=6cm]{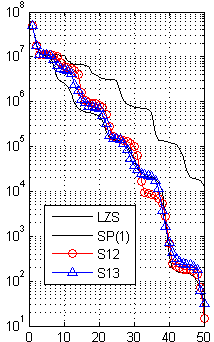}
\includegraphics[width=3.5cm, height=6cm]{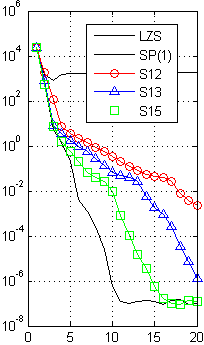}
\caption{Lanczos and SPAM(1) compared with SPAM(1,$\ell$) for small values of $\ell$. Left: reaction-diffusion problem, smallest eigenvalue. Other pictures: largest eigenvalue. Middle left: banded matrix; middle right: bcsstk07; right: bcsstk10. The graphs for Lanczos and SPAM(1) can also be found in Figures \ref{Figure5}, \ref{Figure7}, and \ref{Figure8}.}\label{Figure9}
\end{figure}

\subsection{Comparing SPAM(1,$\ell$) with one-step preconditioned Jacobi-Davidson}

In the $k$-th inner iteration of SPAM(1,$\ell$), the Jacobi-Davidson correction equation (\ref{jdspam}) for $A_k$ is solved using $\ell$ steps of MinRES. We will now compare this with the Jacobi-Davidson method with one-step approximation as preconditioner, as was described in Section \ref{osjd}. This means that in each inner iteration, the initial approximation $A_0$ is used instead of $A_k$. We will still apply $\ell$ steps of MinRES to solve the corresponding correction equation and denote the resulting method by JD(1,$\ell$). Since one of the aims of SPAM was to save on the costs of the matrix-vector products in the inner iteration, we will also apply Jacobi-Davidson {\em without} preconditioning, and approximate the exact correction equation (\ref{aug}) with matrix $A$ by performing $\ell$ steps of MinRES as well, and denote this method by JD($\ell$). Thus, the differences between these three methods lie in the inner iteration: SPAM(1,$\ell$), JD(1,$\ell$) and JD($\ell$) all apply $\ell$ steps of MinRES per inner iteration step, to a linear equation with matrix $A_k$,  $A_0$, and $A$, respectively. Note that in \cite{ShWaTiMi}, no explicit comparison of SPAM with Jacobi-Davidson was made, even though the methods are so closely related.

\begin{figure}[ht]
\includegraphics[width=3.5cm, height=6cm]{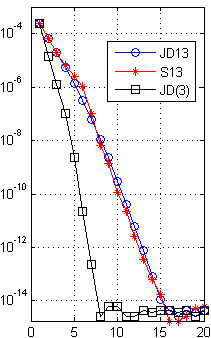}
\includegraphics[width=3.5cm, height=6cm]{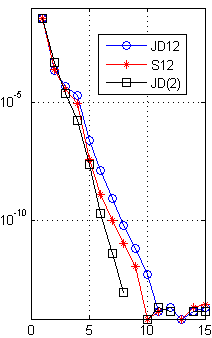}
\includegraphics[width=3.5cm, height=6cm]{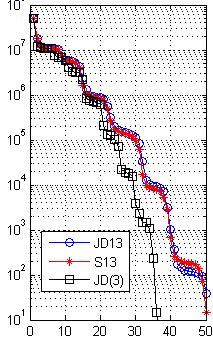}
\includegraphics[width=3.5cm, height=6cm]{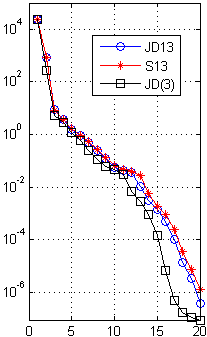}
\caption{Comparing SPAM(1,$\ell$) with JD(1,$\ell$) and JD($\ell$). The eigenvalue problems are exactly the same as the corresponding ones in Figure \ref{Figure9}, and the curves for SPAM(1,$\ell$) can be found there as well.}\label{Figure10}
\end{figure}

As expected, JD($\ell$) is the clear winner in all experiments, although the difference with JD(1,$\ell$) and SPAM(1,$\ell$) is not enough to automatically disqualify the latter two. Since the matrix-vector products in their inner iterations are, in general, considerably cheaper than in JD($\ell$), both methods could be competitive. Having said this, the difference between JD(1,$\ell$) and SPAM(1,$\ell$) is quite small and not always in favor for SPAM(1,$\ell$), even though SPAM(1,$\ell$), much more so than JD(1,$\ell$) uses the information that is available from the outer iteration also in its inner iteration. As argued already in Section \ref{sect2.5.4}, this may actually be less effective than using only the best information that is available from the outer iteration, as the Jacobi-Davidson method does. So far, the numerical experiments are in favor of using Jacobi-Davidson with one-step preconditioning instead of SPAM(1,$\ell$).

\section{Conclusions}
The experiments above illustrate mathematical aspects of SPAM as a method for approximating eigenvalues of a Hermitian matrix. Using approximations from below, SPAM can be seen as a boosted version of the Lanczos method in the sense that convergence towards the largest eigenvalues is stimulated. Since Lanczos itself is often used to provide a good start vector for Jacobi-Davidson, SPAM is a therefore a good candidate for this task, too. Since the difference between SPAM and Jacobi-Davidson with one step approximation is small, it may be preferred to use the latter, especially since the latter is even more easy to use. There does not seem to be a significant gain in re-using the action of $A$ on the orthogonal complement $\UU$ of the current Ritz vector $u$ within $\VV$ also in the inner iterations in comparison with only re-using the action of $A$ on $u$, as Jacobi-Davidson with one step approximation does. This does not mean that the original idea of the authors \cite{ShWaTiMi} of SPAM, to save on the costs of the inner iterations of for instance Jacobi-Davidson, was incorrect. It may well pay off to do so, but this may be done with Jacobi-Davidson with one step approximation just as well. Thus, the main conclusion of this paper is that the value of SPAM probably lies in providing good initial approximations for the Jacobi-Davidson method.

\subsection*{About the dedication}
This paper is dedicated to Gerard Sleijpen's $60$-th birthday. Sleijpen supervised Jan Brandts in both his MSc (1990) and PhD (1995) work, whereas in 2006 he supervised Ricardo Reis da Silva in his MSc work. Both authors recall with pleasure the personal as well as the mathematical interactions.

\end{document}